\documentclass[11pt,reqno]{amsart}
\usepackage{fullpage}
\usepackage{amsmath, amssymb, amsthm,mathrsfs}
\usepackage{enumerate}
\usepackage[foot]{amsaddr}
\usepackage{algorithm}
\usepackage {algpseudocode}
\usepackage{graphicx}
\usepackage{color}
\usepackage{bm}
\usepackage{bigints}
\usepackage{textgreek}
\usepackage[colorlinks=true,linkcolor=red,citecolor=blue,urlcolor=cyan]{hyperref}
\usepackage{verbatim}
\usepackage{enumitem}
\usepackage{mathtools}
\usepackage{longtable}

\theoremstyle{remark}

\usepackage[backend=bibtex,giveninits=true,sorting=nyt,natbib=true,maxcitenames=2,maxbibnames=10,url=false,doi=true,backref=false]{biblatex}

\renewbibmacro{in:}{\ifentrytype{article}{}{\printtext{\bibstring{in}\intitlepunct}}}

\makeatletter
\DeclareFontFamily{U}{tipa}{}
\DeclareFontShape{U}{tipa}{m}{n}{<->tipa10}{}
\newcommand{\arc@char}{{\usefont{U}{tipa}{m}{n}\symbol{62}}}%

\newcommand{\arc}[1]{\mathpalette\arc@arc{#1}}

\newcommand{\arc@arc}[2]{%
  \sbox0{$\m@th#1#2$}%
  \vbox{
    \hbox{\resizebox{\wd0}{\height}{\arc@char}}
    \nointerlineskip
    \box0
  }%
}
\makeatother

\usepackage{pgf}

\begin{document}
\title{Dynamics and stochastic resonance in a mathematical model of bistable phosphorylation and nuclear size control}
\author{Xuesong Bai}
\author{Jonathan Touboul}
\author{Thomas G. Fai}
\email{jtouboul@brandeis.edu, tfai@brandeis.edu}
\address{Department of Mathematics, Brandeis University, Waltham, MA}

\date{\today}

\begin{abstract} 

Robust oscillations play crucial roles in a wide variety of biological processes and are often generated by deterministic mechanisms. However, stochastic fluctuations often generate complex perturbations of these deterministic oscillations, potentially strengthening or weakening their robustness. In this paper, we study bistable phosphorylation as a mechanism for robust oscillation. We present a simple nucleocytoplasmic transport and cell growth model where cargo proteins undergo bistable phosphorylation prior to nuclear import. We perform a detailed bifurcation analysis to examine the system's dynamical behavior. We then introduce additive noise into the model and study the stochastic resonance behavior and robustness of oscillations under noise. Our results show that, depending on the phosphorylation threshold, time-scale parameters, and nucleocytoplasmic transport rate, bistable phosphorylation may generate oscillations via Hopf bifurcations; moreover, stochastic resonance and Bautin bifurcations enhance the robustness of the oscillations.

\end{abstract}

\maketitle

\section{Introduction}

Robust oscillations with precise frequency and amplitude play crucial roles in a wide variety of biological processes, such as circadian rhythm \cite{kim2014molecular,kim2021tick}, cardiac pacemaking \cite{grudzinski2004modeling,maltsev2014modern}, and rhythmic motor behaviors \cite{kintos2008modeling}. Achieving robust oscillation generation is challenging because stochastic fluctuations from noise are ubiquitous in biological systems and span multiple spatial and temporal scales \cite{tsimring2014noise}; therefore, oscillation generation mechanisms must operate under a broad range of perturbations. Identification and characterization of robust mechanisms for oscillation generation are often interesting yet demanding topics in mathematical biology.

Mathematically speaking, oscillations in biology have often been associated with deterministic mechanisms related to the interactions between components \cite{kruse2005oscillations} that are classically modeled by ordinary differential equations (ODEs) \cite{li2018systems}, and can be clearly distinguished from stochastic fluctuations \cite{boaretto2021discriminating}. Although one may expect that biological systems must compensate for stochastic fluctuations to maintain steady frequency and amplitude and to generate robust oscillations, ample work has reported a more complex picture. In some cases, noise may have negative effects on the robustness of oscillations, and noise-resistance mechanisms are present in deterministic oscillators to filter out the unwanted effect of stochastic noise \cite{chae2023spatially,vilar2002mechanisms}. By contrast, in other cases, noise enhances the robustness of the oscillations. For example, it is argued that both intercellular coupling and molecular noise are required for stochastic rhythms in the SCN network \cite{ko2010emergence}. Another example is that small-amplitude noises may enhance the robustness of oscillations in excitable systems by means of self-induced stochastic resonance \cite{muratov2005self}. In this case, the noise-induced oscillations arise in parameter regimes away from the bifurcation thresholds, where the corresponding deterministic dynamics show no limit cycles. As a result, noise expands the parameter regime for oscillations.

The specific problem of interest in this paper is the oscillatory dynamics induced by bistable phosphorylation in nucleocytoplasmic transport. Our initial motivation is to identify the sources of fluctuations in the nuclear-to-cytoplasmic volume ratio (N/C ratio). Although there is ample evidence to show that cells maintain a nearly constant N/C ratio under various growth and osmotic conditions \cite{conklin1914cell,moore2019determination,neumann2007nuclear}, recent experimental observations reported rather significant fluctuations \cite{lemiere2022control}. In that paper, the experimental methods available to track the N/C ratio offer limited temporal resolution, resulting in coarse-grained N/C ratio data that also includes some measurement error. In our previous work \cite{bai2025stochastic}, we studied a simplified stochastic gene translation model and showed that noise in the translation process leads to relative fluctuations in the N/C ratio that are proportional to $1/\sqrt{N}$, where $N$ is a system size variable that determines the order of magnitude of the initial protein numbers. Therefore, given the large number of proteins in cells ($N\sim 10^7$ to $10^8$, based on the order of magnitude given in \cite{rollin2023physical}), stochastic translation would appear unlikely to be a significant source of fluctuations in the N/C ratio.

According to the osmotic balance model of the N/C ratio \cite{deviri2022balance,lemiere2022control}, a sufficiently large number of biomolecules such as proteins must enter or exit the nucleus simultaneously to generate significant fluctuations in the N/C ratio. One possible mechanism for proteins to enter the nucleus in large groups is bistable phosphorylation. Some species of nuclear cargo proteins must be phosphorylated prior to being imported into the nucleus, such as the PER protein \cite{beesley2020wake,chae2023spatially}. If the phosphorylation response curve is bistable, cargo proteins may be phosphorylated in a coordinated manner, allowing a large number of them to be phosphorylated and imported simultaneously \cite{chae2023spatially}, potentially leading to oscillations in protein concentrations and the N/C ratio.

Our model in this paper combines the bistable phosphorylation model presented in \cite{chae2023spatially} with the simplified gene translation model from our previous work \cite{bai2025stochastic}. We use an alternative exponential cell-growth model to further simplify the system, so that the ODE system we study contains only two variables. The 2D model allows us to examine the phase plane and the nullclines to further illustrate the mechanisms behind the observed dynamical behaviors.

We perform a detailed bifurcation analysis with respect to the parameters in the phosphorylation and nucleocytoplasmic transport processes to determine the parameter ranges that lead to oscillations and to explore the model's possible dynamical behaviors. Motivated by the original study \cite{chae2023spatially} on robust circadian rhythms, we also consider whether bistable phosphorylation enhances the robustness of oscillations in our simplified model. We introduce noise into the ODE system and study its effect on oscillations. In particular, we show that bistable phosphorylation can lead to self-induced stochastic resonance \cite{muratov2005self}, allowing the system to oscillate over a wider parameter range in the presence of noise.

\section{Mathematical model of bistable phosphorylation in nuclear import and cell growth}

The model we study in this paper is a combination of a cell growth model from our previous work \cite{bai2025stochastic} and a bistable phosphorylation model based on a circadian rhythm model described in~\cite{chae2023spatially}. The model is based on the following two considerations:
\begin{enumerate}
\item In \cite{chae2023spatially}, there are two potential rhythm generation mechanisms: bistable phosphorylation of PER proteins and negative feedback on the translation rate of PER proteins. In this paper, we focus on the standalone effect of bistable phosphorylation on oscillation generation; therefore, we omit negative feedback mechanisms from the present model.
\item The original circadian rhythm model in \cite{chae2023spatially} describes a static cell without volume growth. Building on  our previous studies on the N/C ratio and cell growth \cite{BAI2025112250,bai2025stochastic}, we study dynamical behavior in a growing cell rather than a static one.
\end{enumerate}

\subsection{A simplified model for exponential growth}

In this paper, we use a gene translation model from our previous work \cite{bai2025stochastic} as the basis of the cell growth model. To facilitate the bifurcation analysis in later sections, we formulate a simplified exponential-growth model with fewer variables. Instead of explicitly modeling the autocatalytic replication of ribosomes as in the original gene translation model \cite{bai2025stochastic}, we assume that the translation rate of nuclear cargo proteins and cytoplasmic proteins is proportional to the volume of the cytoplasm:

\begin{subequations}\label{translation2}
\begin{align}
 \label{Pno2}\emptyset&\xrightarrow{k_{vt}\phi_n V_{cyto}}P_{no}, \\
 \label{Pc2}\emptyset&\xrightarrow{k_{vt}\phi_{cy} V_{cyto}}P_{cy},
\end{align}
\end{subequations}
in which $P_{no}$ denotes the nuclear cargo proteins in the cytoplasm, $P_{cy}$ denotes other cytoplasmic proteins, $k_{vt}$ is the translation rate per volume, $\phi_n$ and $\phi_{cy}$ are the gene fractions of $P_{no}$ and $P_{cy}$, respectively, and $V_{cyto}$ is the volume of the cytoplasm, which is positively correlated with the number of $P_{no}$ and $P_c$.

\subsection{Volume growth and number of proteins vs. concentration of proteins}\label{Concentration}

\subsubsection{Numbers vs. concentrations: a note on the model variables}
Throughout this paper, we use $p_j$ and $c_j$ to denote the number and the corresponding concentration of each protein species, respectively. That is, for cytoplasmic species $P_{cy}$, $P_{no}$, and phosphorylated cargo proteins $P_{nop}$ (described in Section \ref{Import}),
\begin{equation*}
 c_j = \frac{p_j}{V_{cyto}}\quad\mbox{for}\quad j = no, nop, cy;
\end{equation*}
and for cargo proteins $P_{ni}$ that have been imported into the nucleus,
\begin{equation*}
 c_j = \frac{p_j}{V_n}\,\mbox{for}\,j = ni,
\end{equation*}
where $V_{cyto}$ is the volume of the cytoplasm mentioned above and $V_n$ is the volume of the nucleus.

\subsubsection{Quotient rule and volume growth}
For simplicity, we assume that the volume of each compartment (i.e., nucleus and cytoplasm) is proportional to the number of proteins in that compartment:
\begin{equation}\label{cellvolume}
V_{cyto} = A_{cyto}(p_{cy}+p_{no}+p_{nop}),
\end{equation}
and
\begin{equation}
V_n = A_n p_{ni},
\end{equation}
where $A_{cyto}$ and $A_n$ are scaling factors for the cytoplasmic and nuclear volume, respectively. This proportionality is motivated by our previous study of volume control by osmotic forces in \emph{S. pombe} \cite{lemiere2022control}, under the simplifying assumption that $P_{ni}$ is the dominant nuclear osmolyte. Since the growth of the number of proteins in (\ref{translation2}) is exponential, the number of proteins and volumes $V_{cyto}$, $V_n$ will approach infinity as $t\to\infty$, and the system will not have any finite steady state, making the ODE system unsuitable for bifurcation analysis. Therefore, we choose to use protein concentrations as variables in our ODE system. Later in Section \ref{results}, we show  that the ODE system of protein concentrations either converges to a finite set of steady states or exhibits oscillations, depending on the parameters.

We then use the quotient rule to calculate the time derivatives of the concentrations:
\begin{equation}
\frac{dc_j}{dt} = \frac{d}{dt}\left(\frac{p_j}{V_k}\right) = \frac{V_k\frac{dp_j}{dt}-p_j\frac{dV_k}{dt}}{V_k^2}=\frac{1}{V_k}\frac{dp_j}{dt}-\frac{p_j}{V_k^2}\frac{dV_k}{dt},
\end{equation}
where the index is over $k=cyto$ for $j=cy,no,nop$ and $k=n$ for $j=ni$. The ODEs for the protein numbers (i.e., $dp_j/dt$) follow from the law of mass action (\ref{translation2}), and from the propensity functions specified later in (\ref{phosphorylation}).

The N/C ratio $\Phi_{NCyto}$ is defined by $\Phi_{NCyto} = V_n/V_{cyto}$, and we also use the quotient rule to calculate $d\Phi_{NCyto}/dt$.

\subsection{Mathematical model of bistable phosphorylation in nuclear import}\label{Import}
We assume that nuclear cargo proteins must be phosphorylated before they are imported into the nucleus:
\begin{subequations}\label{phosphorylation}
\begin{align}
P_{no}&\xrightarrow{f_p(c_{no},c_{nop},K)}P_{nop}, \\
P_{nop}&\xrightarrow{f_{dp}(c_{no},c_{nop})}P_{no}, \\
P_{nop}&\xrightarrow{k_{nt}}P_{ni}.
\end{align}
\end{subequations}
Here, $P_{no}$ denotes the non-phosphorylated nuclear cargo proteins in the cytoplasm, $P_{nop}$ denotes the phosphorylated nuclear cargo proteins in the cytoplasm, and $P_{ni}$ denotes the nuclear cargo proteins that have been imported into the nucleus. Note that the variables $c_{no}$ and $c_{nop}$ in the above propensity functions are concentrations of the corresponding proteins $P_{no}$ and $P_{nop}$, respectively, while the original version of the functions in \cite{chae2023spatially} uses numbers of proteins as variables. 

We use the propensity functions specified in \cite{chae2023spatially} to model bistable phosphorylation:
\begin{equation}\label{fp}
 f_p(c_{no},c_{nop},K_c,\tau) = \frac{1}{\tau}\left[\frac{(c_{no}+c_{nop})^{m+1}}{(c_{no}+c_{nop})^m+[f_{sca}(c_{nop},K_c)]^m}\right],
\end{equation}
\begin{equation}\label{fdp}
 f_{dp}(c_{no},c_{nop},\tau) = \frac{1}{\tau}c_{nop},
\end{equation}
and the scaling function $f_{sca}$ is given by
\begin{equation}\label{fsca1}
 f_{sca}(c_{nop},K_c) = K_c\max\left\{\frac{1}{8},\frac{5}{4} - \frac{4c_{nop}}{K_c}\right\}.
\end{equation}
Note that (\ref{fsca1}) is not differentiable at $c_{nop}=9K_c/32$. This non-differentiability may adversely affect the numerical bifurcation analysis. Moreover, smoothness of ODE systems is often required in theoretical studies of dynamical systems \cite{kuznetsov1998elements}.  For these two reasons, we come up with a similar-shaped but differentiable scaling function in the numerical simulations:
\begin{equation}\label{fsca2}
 f_{sca}(c_{nop},K_c) = K_c\left[\frac{5}{4}-\frac{9}{8}\frac{(c_{nop}/K_c)^{m_{sca}}}{(c_{nop}/K_c)^{m_{sca}}+(9/64)^{m_{sca}}}\right].
\end{equation}

\subsection{ODE system used for bifurcation analysis}

Reactions (\ref{translation2}) and (\ref{phosphorylation}), together with the nuclear and cytoplasmic volume model (\ref{cellvolume}), comprise the complete bistable phosphorylation and cell growth model. The corresponding ODE system (\ref{odeexp2}) contains seven ODEs:
\begin{subequations}\label{odeexp2}
 \begin{align}
  \label{ode1exp2} \frac{dc_{no}}{dt} &= k_{vn}-f_p(c_{no},c_{nop},K_c,\tau)+f_{dp}(c_{no},c_{nop},\tau)-c_{no}A_{cyto}(k_{vn}+k_{vcy}-k_{nt}c_{nop}), \\
  \label{ode2exp2} \frac{dc_{nop}}{dt} &= f_p(c_{no},c_{nop},K_c,\tau)-f_{dp}(c_{no},c_{nop},\tau)-k_{nt}c_{nop}-c_{nop}A_{cyto}(k_{vn}+k_{vcy}-k_{nt}c_{nop}), \\
  \label{ode3exp2} \frac{dc_{ni}}{dt} &= \frac{1}{\Phi_{NCyto}}k_{nt}c_{nop}(1-A_nc_{ni}),\\
  \label{ode4exp2} \frac{dc_{cy}}{dt} &= k_{vcy}-c_{cy}A_{cyto}(k_{vn}+k_{vcy}-k_{nt}c_{nop}), \\
  \label{ode5exp2} \frac{d\Phi_{NCyto}}{dt} &= A_nk_{nt}c_{nop}-A_{cyto}(k_{vn}+k_{vcy}-k_{nt}c_{nop})\Phi_{NCyto}, \\
  \label{ode6exp2} \frac{dV_n}{dt} &= A_nk_{nt}c_{nop}V_{cyto}, \\
  \label{ode7exp2} \frac{dV_{cyto}}{dt} &= A_{cyto}(k_{vn}+k_{vcy}-k_{nt}c_{nop})V_{cyto},
 \end{align}
\end{subequations}
where $k_{vn} = k_{vt}\phi_n$ and $k_{vcy} = k_{vt}\phi_{cy}$. Note that if we assume that the gene fractions satisfy $\phi_n+\phi_{cy}=1$, then $k_{vn}+k_{vcy}=k_{vt}$.

For the purpose of bifurcation analysis, we further reduce the system to focus on rhythm-generation mechanisms. We note that the only potential mechanism in the system that generates oscillations is the bistable phosphorylation that appears in equations (\ref{ode1exp2}) and (\ref{ode2exp2}), and these two ODEs form a self-contained subsystem that depends only on $c_{no}$ and $c_{nop}$. Based on the observations above, to perform a bifurcation analysis of the system and study its oscillatory behaviors, it suffices to focus on the ODEs that involve the synthesis of nuclear cargo proteins and bistable phosphorylation. That is, we study the bifurcation of the 2D ODE system consisting of (\ref{ode1exp2}) and (\ref{ode2exp2}).

\section{Dynamical behaviors and bifurcation analysis}\label{results}

The codes for the following results are available on GitHub\footnote{\href{https://github.com/topgunbai683/bistable-phosphorylation-dynamics.git}{https://github.com/topgunbai683/bistable-phosphorylation-dynamics.git}}.

\subsection{Fraction of phosphorylated cargo proteins and analysis of dynamical behaviors}

In this section, we present examples of trajectories and nullclines in the phase plane to illustrate the mechanisms behind the various dynamical behaviors and bifurcations.

Instead of using the default $c_{no}$ - $c_{nop}$ plane, we plot the fraction $c_{nop}/(c_{no}+c_{nop})$ of phosphorylated cargo proteins against the total concentration $c_{no}+c_{nop}$ of cargo proteins. The motivation is that the phosphorylation function (\ref{fp}) is controlled by $c_{no}+c_{nop}$. If the timescale parameter $\tau$ is small, then the fraction $c_{nop}/(c_{no}+c_{nop})$ quickly equilibrates to
\begin{equation*}
 \frac{c_{nop}}{c_{no}+c_{nop}} \approx \frac{(c_{no}+c_{nop})^m}{(c_{no}+c_{nop})^m+[f_{sca}(c_{nop},K_c)]^m},
\end{equation*}
which is an S-shaped response curve if the scaling function $f_{sca}$ is set to be equal to either of the expressions defined in (\ref{fsca1}) or (\ref{fsca2}) \cite{chae2023spatially}. Therefore, plots in the $(c_{no}+c_{nop})$ - $c_{nop}/(c_{no}+c_{nop})$ plane clearly show the bistability in the phosphorylation process, and demonstrate how variations in the parameters affect the nullclines and change the dynamical behaviors of the system (see Fig. \ref{KcExamples}(a), for example). 

\subsubsection{Role of the phosphorylation threshold \texorpdfstring{$K_c$}{Kc}}
\begin{figure}[h]
 \includegraphics[width=0.9\textwidth]{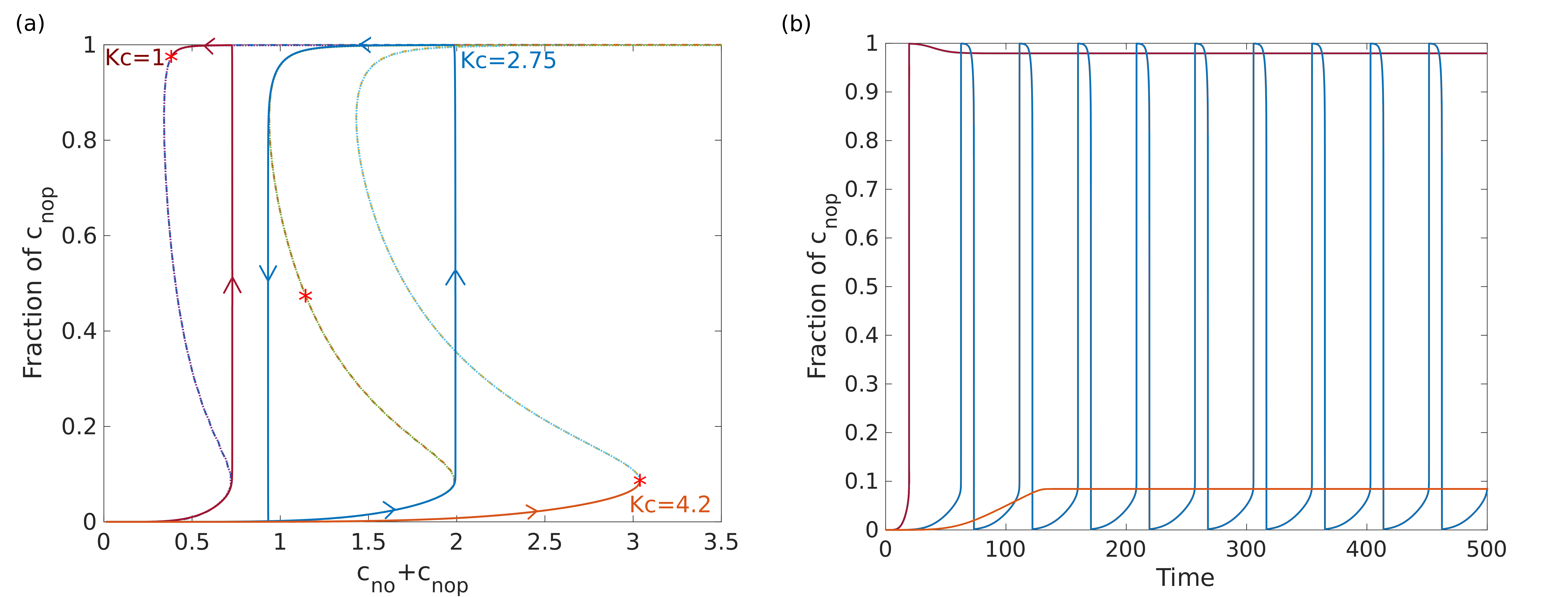}
 \caption{Effect of $K_c$ on the dynamical behaviors of the system. \textbf{(a)}, phase plane. Solid lines, trajectories; dashed lines, $c_{no}$-nullclines; dotted lines, $c_{nop}$-nullclines. For three values of $K_c$: the nullclines with leftmost inflection point and bordeaux trajectories correspond to $K_c=1$, the nullclines with the middle inflection point and blue periodic trajectory correspond to $K_c=2.75$, while the nullclines with the rightmost inflection point and red trajectory correspond to $K_c=4.2$. In each case, there was a single fixed point, marked with a red asterisk. \textbf{(b)}, corresponding trajectories of the $c_{no}$ fraction as a function of time. The colors match those used in panel (a): bordeaux for $K_c=1$, blue for $K_c=2.75$ and red for $K_C=4.2$. In all simulations, $\tau=0.01$ and $k_{nt}=0.1$. }
 \label{KcExamples}
\end{figure}

Fig. \ref{KcExamples} shows three distinct types of dynamical behavior of the system, as the phosphorylation threshold $K_c$ varies. Note that both the $c_{no}$- and the $c_{nop}$-nullclines remain S-shaped and bistable in the examples, with $\tau=0.01$ and $k_{nt}=0.1$. As $K_c$ increases, the widths of the S-shaped regions of the nullclines increase, but the overall shapes of the nullclines are unaffected.

Although the two nullclines are nearly indistinguishable in this representation, they have a single intersection whose location varies with Kc (red asterisks in Figure~\ref{KcExamples}). We summarize the possible cases of intersection and the resulting dynamical behaviors in the following.

\begin{enumerate}
\item When $K_c$ is small (Fig. \ref{KcExamples}, $K_c=1$), the two nullclines intersect at the top stable branches of the S-shapes, resulting in one stable equilibrium with a high fraction of phosphorylated cargo proteins (bottom region in Fig. \ref{tauKcBifurcation}).
\item When $K_c$ is large (Fig. \ref{KcExamples}, $K_c=4.2$), the two nullclines intersect at the lower stable branches of the S-shapes, resulting in one stable equilibrium with a low fraction of phosphorylated cargo proteins (top region in Fig. \ref{tauKcBifurcation}).
\item In between the two cases (Fig. \ref{KcExamples}, $K_c=2.75$), the two nullclines intersect at the unstable branches of the S-shapes, the only equilibrium is thus unstable, and the trajectories converge to a stable limit cycle.
\end{enumerate}

\subsubsection{Role of the phosphorylation timescale parameter \texorpdfstring{$\tau$}{tau}}\label{tauanalysis}
\begin{figure}[h]
 \includegraphics[width=0.9\textwidth]{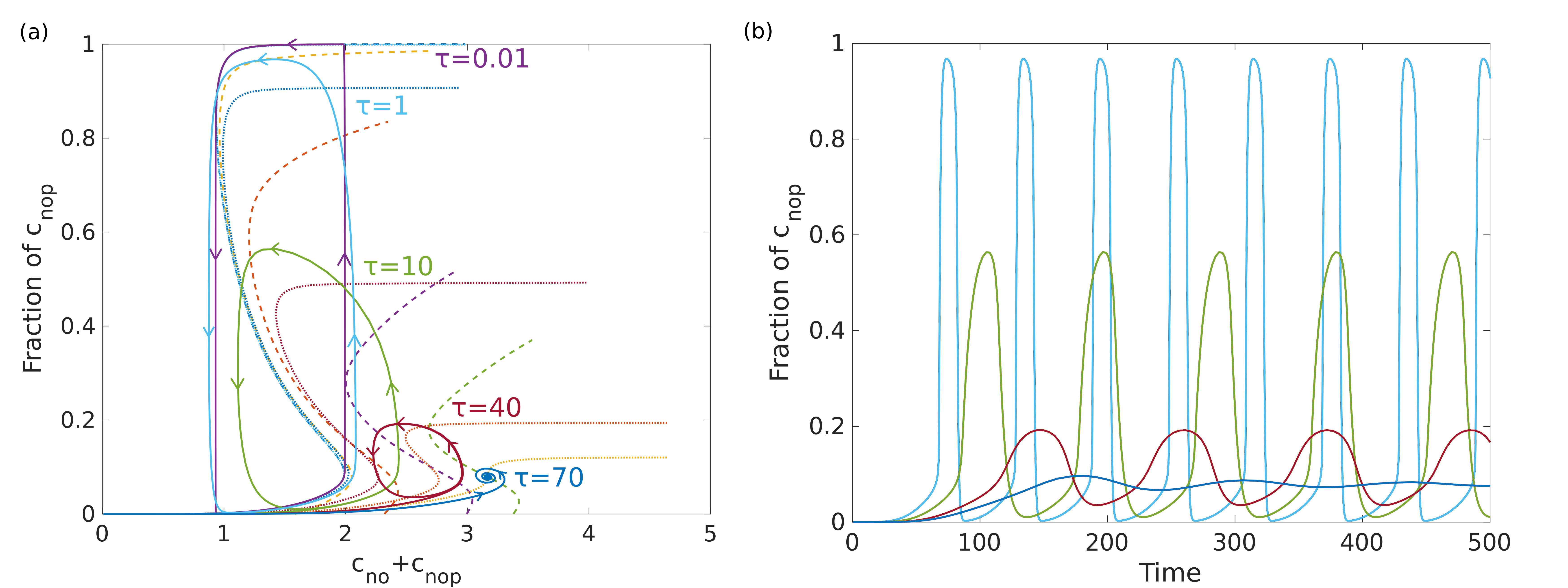}
 \caption{Effect of $\tau$ on the dynamical behaviors of the system. \textbf{(a)}, phase plane. Solid lines, trajectories; dashed lines, $c_{no}$-nullclines; dotted lines, $c_{nop}$-nullclines. For five values of $\tau$: the purple trajectory corresponds to $\tau=0.01$, the light blue trajectory corresponds to $\tau=1$, the green trajectory corresponds to $\tau=10$, the bordeaux trajectory corresponds to $\tau=40$, and the blue trajectory corresponds to $\tau=70$. \textbf{(b)}, corresponding trajectories of $c_{no}$ fraction as a function of time. The colors match those used in panel (a): light blue for $\tau=1$, green for $\tau=10$, and bordeaux for $\tau=40$. In all simulations, $K_c=2.75$ and $k_{nt}=0.1$.}
 \label{tauExamples}
\end{figure}

Timescales are particularly important for the generation of robust oscillations. Fig. \ref{tauExamples} illustrates how varying the phosphorylation timescale parameter $\tau$ alters the dynamics and trajectories of the system. 

For very small $\tau$ (Fig. \ref{tauExamples}, $\tau=0.01$), the right-hand sides of the ODE system (\ref{ode1exp2}) and (\ref{ode2exp2}) are dominated by the phosphorylation function $f_p$ and dephosphorylation function $f_{dp}$. As a result, the $c_{no}$- and $c_{nop}$- nullclines are very close to each other and look nearly identical. Moreover, due to rapid phosphorylation and dephosphorylation, when the trajectory is traveling along one stable branch and hits the corresponding fold point of the S-shaped nullcline, it will instantly equilibrate to the opposite stable branch before any significant change in $c_{no}+c_{nop}$ occurs. Consequently, the trajectory features near-vertical segments and sharp corners.

As $\tau$ increases, the phosphorylation functions become less dominant, and the $c_{no}$- and $c_{nop}$- nullclines become more distinct. The trajectory approaches the opposite stable branch more slowly at the fold points, and the time for this process is long enough for $c_{no}+c_{nop}$ to have noticeable variations. As a result, the trajectory no longer has near-vertical segments and sharp corners (Fig. \ref{tauExamples}, $\tau=1$).

If $\tau$ increases further, the phosphorylation process gradually slows down enough to become a limiting factor in the dynamics. The fraction of $c_{nop}$ becomes smaller and smaller, so that the maximum value of the $c_{nop}$- nullcline decreases, leading to lower amplitude in the oscillation (Fig. \ref{tauExamples}, $\tau=10$ and $\tau=40$). Eventually, the $c_{nop}$- nullcline is no longer S-shaped and becomes step-like. The system no longer has oscillations in this case; instead, there exists one stable equilibrium (Fig. \ref{tauExamples}, $\tau=70$).

\subsubsection{Role of the nucleocytoplasmic transport rate \texorpdfstring{$k_{nt}$}{knt}}
\begin{figure}[h]
 \includegraphics[width=0.45\textwidth]{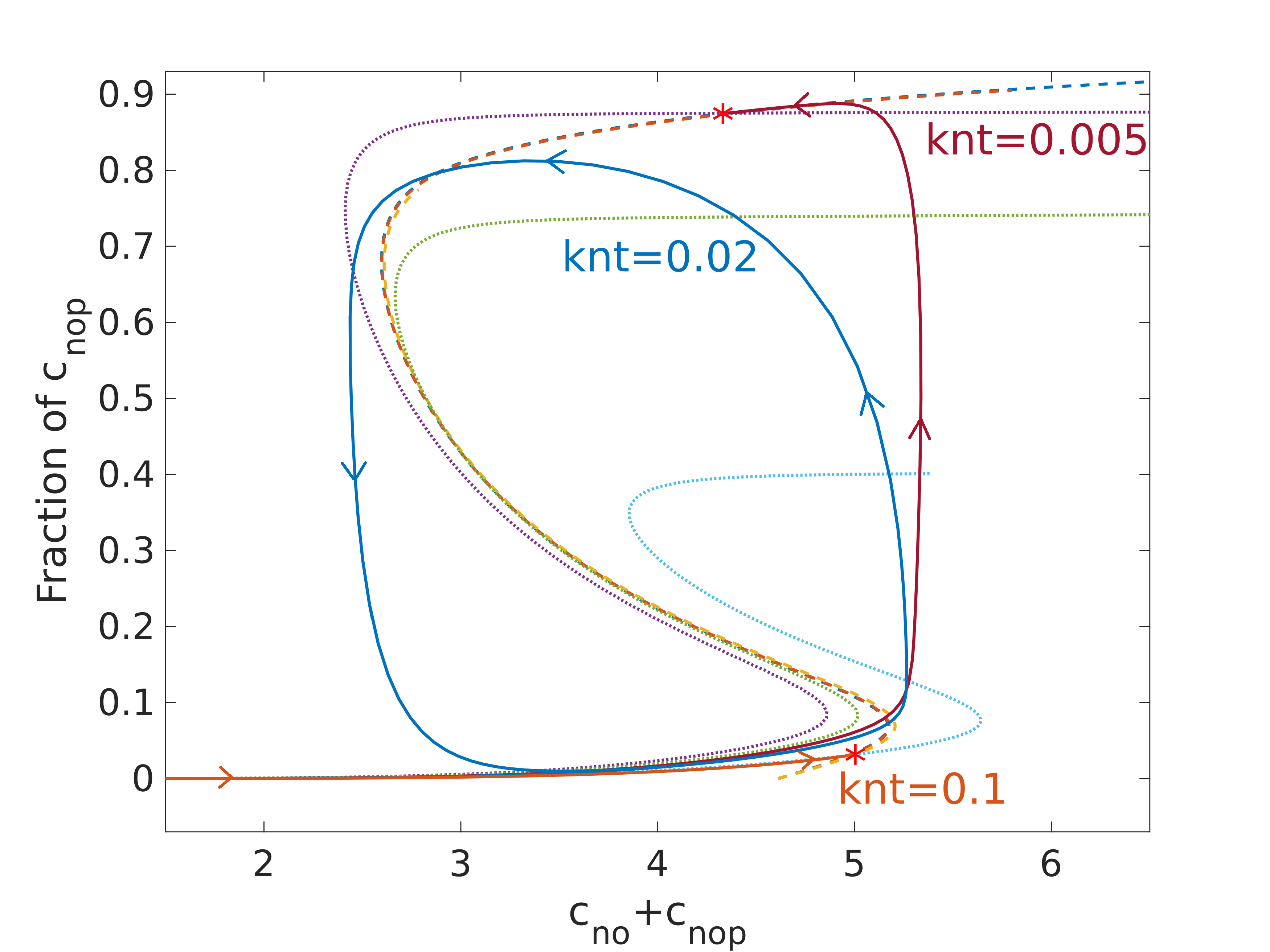}
 \caption{Effect of $k_{nt}$ on the dynamical behaviors of the system. Solid lines, trajectories; dashed lines, $c_{no}$-nullclines; dotted lines, $c_{nop}$-nullclines. For three values of $k_{nt}$: the $c_{nop}$-nullcline with leftmost inflection point and bordeaux trajectories correspond to $k_{nt}=0.005$, the $c_{nop}$-nullcline with the middle inflection point and blue periodic trajectory correspond to $k_{nt}=0.02$, while the $c_{nop}$-nullcline with the rightmost inflection point and red trajectory correspond to $k_{nt}=0.1$. In each case, there was a single fixed point, marked with a red asterisk. In all simulations, $K_c=6.56$ and $\tau=14.5$.}
 \label{kntExamples}
\end{figure}

Fig. \ref{kntExamples} shows three distinct types of dynamical behavior of the system, as the nucleocytoplasmic transport rate $k_{nt}$ varies. Again, both the $c_{no}$- and the $c_{nop}$-nullclines remain S-shaped and bistable in the examples, with $K_c=6.56$ and $\tau=14.5$. Note that $k_{nt}$ has little effect on the $c_{no}$- nullcline; we see that the three $c_{no}$- nullclines in Fig. \ref{kntExamples} are almost indistinguishable. However, $k_{nt}$ has significant effects on the $c_{nop}$- nullcline.

As $k_{nt}$ increases, phosphorylated cargo proteins are imported more rapidly into the nucleus, and the fraction of $c_{nop}$ decreases. Consequently, the maximum value of the $c_{nop}$ nullcline decreases, altering the system's behavior.

\begin{enumerate}
\item When $k_{nt}$ is small (Fig. \ref{kntExamples}, $k_{nt}=0.005$), the two nullclines intersect at the top stable branches of the S-shapes, resulting in one stable equilibrium with a high fraction of phosphorylated cargo proteins (bottom region in Fig. \ref{KckntBifurcation} (a)).
\item When $k_{nt}$ is large (Fig. \ref{kntExamples}, $k_{nt}=0.1$), the two nullclines intersect at the lower stable branches of the S-shapes, resulting in one stable equilibrium with a low fraction of phosphorylated cargo proteins (top region in Fig. \ref{KckntBifurcation} (a)).
\item In between the two cases (Fig. \ref{kntExamples}, $k_{nt}=0.02$), the two nullclines intersect at the unstable branches of the S-shapes, resulting in one unstable equilibrium and one stable limit cycle.
\end{enumerate}

\subsubsection{Bistability of equilibria}\label{BistabilityinEquilibria}
\begin{figure}[h]
 \includegraphics[width=0.45\textwidth]{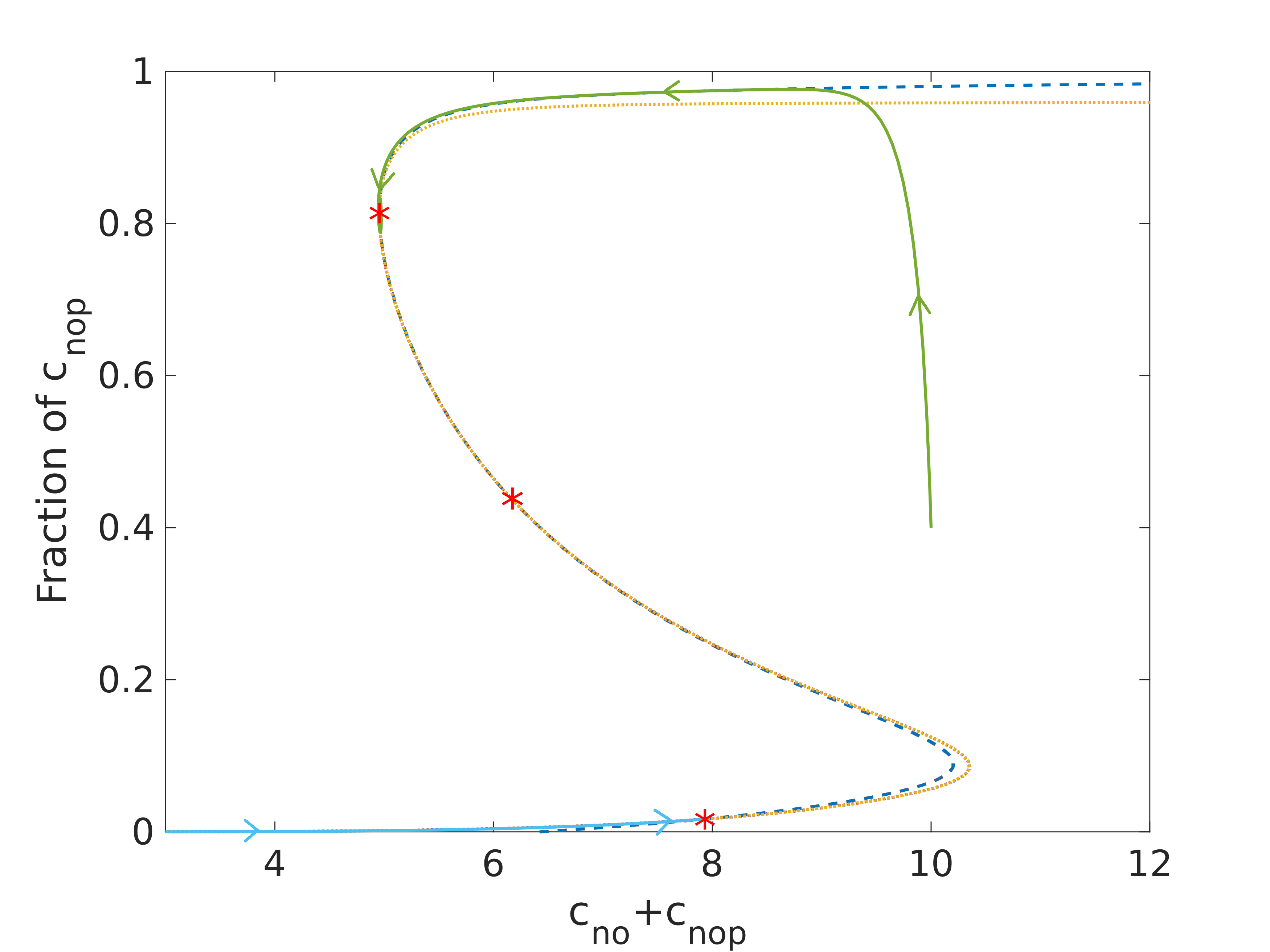}
 \caption{Example of bistability of equilibria. Solid lines, trajectories; dashed lines, $c_{no}$-nullclines; dotted lines, $c_{nop}$-nullclines. Each fixed point is marked with a red asterisk. In all simulations, $K_c=14.2$, $\tau=5$, and $k_{nt}=0.00397$.}
 \label{BistabilityExamples}
\end{figure}

Fig. \ref{BistabilityExamples} shows how the bistability of equilibria occurs due to the positions of the nullclines. We see that there are three points of intersection between the two nullclines. Two of the intersections occur at the stable branches, one at the top and the other at the bottom, resulting in two stable equilibria. The other intersection occurs at the unstable branch, leading to an unstable equilibrium.

\subsection{Two-parameter bifurcations: 2D cross sections of the 3D parameter space}

To account for the behaviors of the system illustrated in the previous sections, we now turn to analyzing the bifurcations of the system with respect to the parameters. We examine two-parameter cross-sections of the three parameter space $(\tau,K_c,k_{nt})$, to finely delineate the boundaries between distinct dynamical behaviors and identify various types of oscillatory behavior. For the sake of brevity, we focus our attention particularly on the two-parameter bifurcation diagrams of the system as a function of $(\tau, K_c)$ and $(K_c,k_{nt})$

\subsubsection{Two-parameter bifurcation analysis with respect to \texorpdfstring{$(\tau,K_c)$}{(tau,Kc)}}\label{2parametertauKc}\label{tauKcBifurcationSec}
\begin{figure}[!htbp]
 \includegraphics[width=0.9\textwidth]{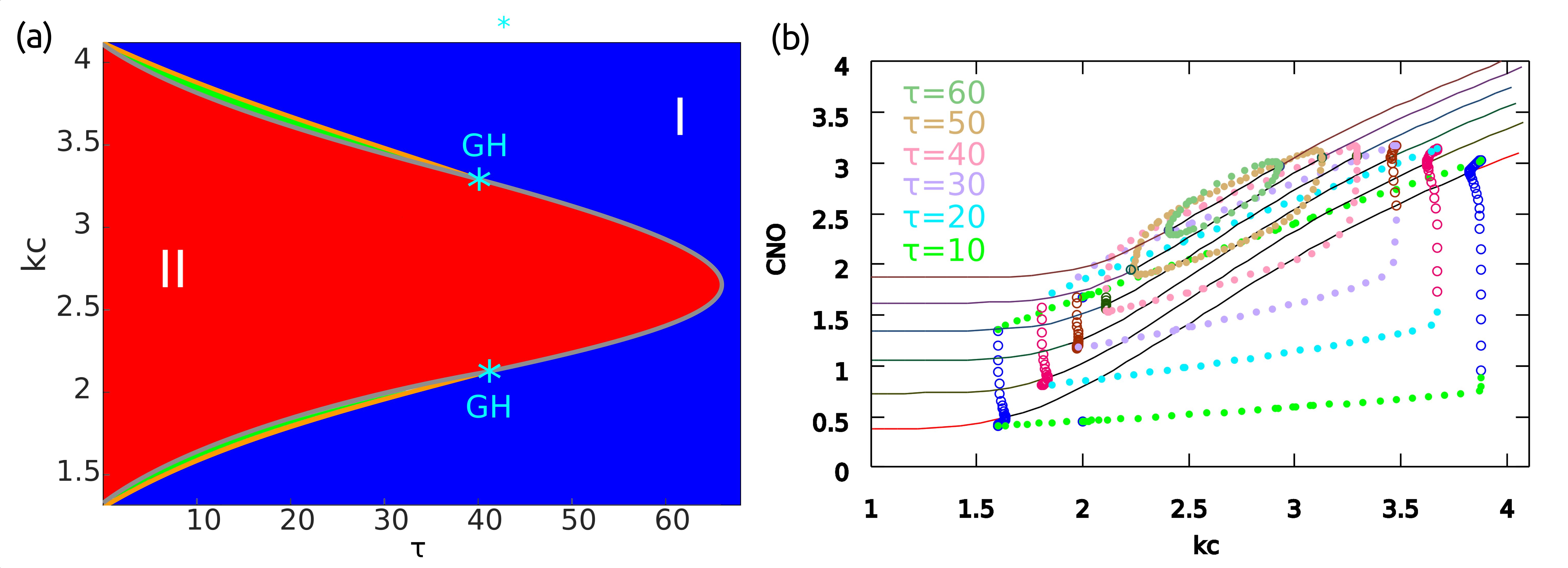}
 \caption{\small Two-parameter bifurcation analysis with respect to $(\tau,K_c)$ with $k_{nt}=0.1$. \textbf{(a)}, $(\tau, K_c)$ plane partitioned by Hopf bifurcation curve of equilibria (gray) and fold bifurcation curve of limit cycles (orange). Blue (Region I), one stable equilibrium. Red (Region II), one stable limit cycle. Green (Region III), two limit cycles of opposite stability. GH, Bautin (generalized Hopf) bifurcation points (labeled by asterisks). Oscillations occur in red and green regions. Diagram computed using the Matcont package on Matlab \cite{dhooge2003matcont}. \textbf{(b)}, bifurcation with respect to $K_c$ under various $\tau$ values. Stable and unstable limit cycles are shown by filled and hollow circles, respectively. Diagram computed using XPPAUT.}
 \label{tauKcBifurcation}
\end{figure}

We first perform a detailed two-parameter bifurcation analysis with respect to the phosphorylation timescale parameter $\tau$ and the phosphorylation threshold $K_c$, where the nucleocytoplasmic transport rate is fixed at the nominal value $k_{nt}=0.1$. The bifurcation diagram featured Hopf bifurcations of equilibria, that switched from sub-critical to supercritical through two Bautin (Generalized Hopf, or GH) bifurcations. The universal unfolding of these codimension two bifurcations includes fold bifurcation curves of limit cycles, which we also tracked numerically. The resulting bifurcation diagram is shown in Fig. \ref{tauKcBifurcation}(a), with a bell-shaped Hopf bifurcation curve depicted in gray and the two Bautin bifurcation points located close to $\tau=40$. The fold bifurcation curves of limit cycles are depicted in orange. As expected from the universal unfolding of the Bautin bifurcation \cite{guckenheimer2013nonlinear,kuznetsov1998elements}, fold bifurcation curves of limit cycles emerge from the Bautin point tangentially to the Hopf bifurcation curve and delineate a region of bistability between oscillations and fixed points (green region in Figure \ref{tauKcBifurcation}(a)) The two fold bifurcation curves of limit cycles then approach the Hopf bifurcation curve again as $\tau\to 0$. 

The bifurcation curves partition the 2D $(\tau, K_c)$ parameter space into regions with distinct dynamical behaviors. To the right of the Bautin points, the Hopf bifurcations are supercritical: they have negative first Lyapunov coefficients and are associated with stable limit cycles (Fig. \ref{tauKcBifurcation} (b), $\tau=50$ and $\tau=60$); whereas to the left of the Bautin points, the Hopf bifurcation is subcritical, associated with a positive first Lyapunov coefficient, and generates unstable limit cycles (Fig. \ref{tauKcBifurcation} (b), $\tau=10$ to $\tau=30$); each unstable limit cycle coexists with a corresponding stable limit cycle in Region III of Fig. \ref{tauKcBifurcation}(a), until crossing the fold bifurcation curve of limit cycles and entering Region I where the system has unique stable equilibria and no oscillations, or until crossing the Hopf bifurcation curve of equilibria and entering Region II where where the system has unique stable limit cycles.

Since Region III (coexistence of stable and unstable limit cycle solutions) is narrow compared to Region II (unique stable limit cycle solutions), the sizes of the limit cycles increase very rapidly in Region III, leading to an apparent sudden onset of oscillations with nonzero amplitude when the phosphorylation threshold $K_c$ crosses the bifurcation values (Fig. \ref{tauKcBifurcation}(b), $\tau=10$ to $\tau=30$). 

\begin{figure}[h]
 \includegraphics[width=0.6\textwidth]{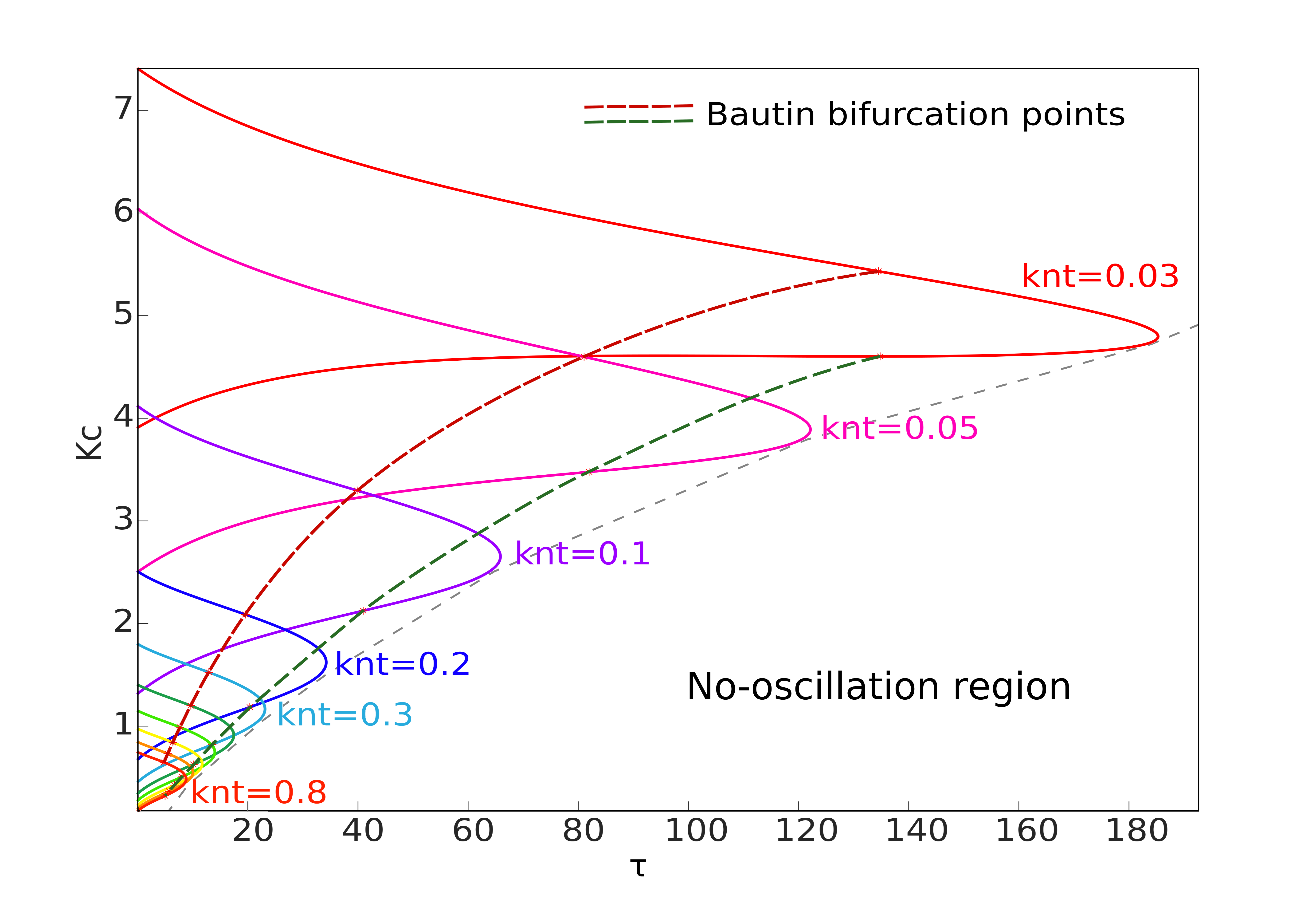}
 \caption{Hopf bifurcation curves of equilibria in the $(\tau,K_c)$ plane with $k_{nt}=$ 0.03, 0.05, 0.1, 0.2, 0.3, 0.4, 0.5, 0.6, 0.7, and 0.8. The top and bottom Bautin bifurcation points on each curve are connected by dark red and dark green dashed curves, respectively.}
 \label{tauKcHopfCurves}
\end{figure}

Next, we track the location of the two-parameter Hopf bifurcations curves as a function of the nucleocytoplasmic transport rate $k_{nt}$. Fig. \ref{tauKcHopfCurves} represents the two-parameter Hopf bifurcation curves and the Bautin bifurcations (marked by dashed curves) for a variety of $k_{nt}$ values to identify the structure of the three-parameter Hopf bifurcation surface in the 3D parameter space through its cross-sections parallel to the $(\tau, K_c)$ plane.

As $k_{nt}$ increases, the region enclosed by the Hopf bifurcation curve (corresponding to Region II in Fig. \ref{tauKcBifurcation} (a)) shrinks. That is, the region of oscillations gets smaller with larger $k_{nt}$; in this case, both $\tau$ and $K_c$ need to decrease for oscillations to occur. For small $k_{nt}$ ($\sim 10^{-2}$), $\tau$ can be quite large ($\sim 10^2$) while oscillations still persist; whereas for large $k_{nt}$ ($\sim 10^{-1}$), $\tau$ is restricted to much smaller values ($\sim 10^1$) to maintain oscillations. Moreover, there exists a distinct region in the $(\tau, K_c)$ plane where oscillations are impossible, no matter the value of $k_{nt}$; these are associated with regimes of slow phosphorylation timescale and low threshold. We analyze the effect of large $\tau$ and the mechanism behind this no-oscillation region in Section \ref{tauanalysis} above.

\subsubsection{Two-parameter bifurcation analysis with respect to \texorpdfstring{$(K_c,k_{nt})$}{(Kc,knt)}}\label{2parameterKcknt}
\begin{figure}[h]
 \includegraphics[width=0.9\textwidth]{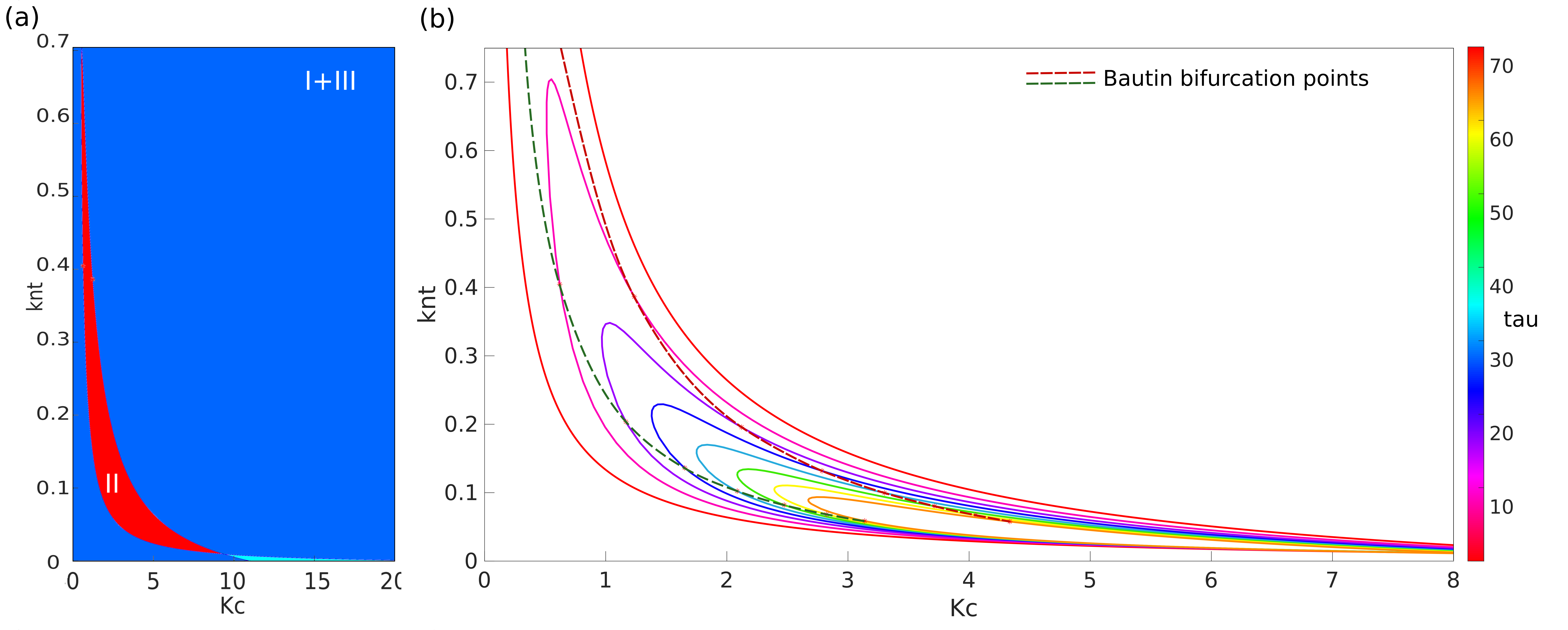}
 \caption{Two-parameter bifurcation analysis with respect to $(K_c,k_{nt})$. \textbf{(a)}, $(K_c,k_{nt})$ plane partitioned by Hopf bifurcation curve of equilibria, with $\tau=10$. Blue (Region I+III), one stable equilibrium except for very narrow regions along the Hopf bifurcation curve, where stable and unstable limit cycles coexist. Red (Region II), one stable limit cycle. Cyan (Region IV), two stable equilibria and one unstable equilibrium (bistability). \textbf{(b)}, Hopf bifurcation curves of equilibria in the $(K_c,k_{nt})$ plane from $\tau=0.01$ to $\tau=70$. The top and bottom Bautin bifurcation points on each curve are connected by dark red and dark green dashed curves, respectively.}
 \label{KckntBifurcation}
\end{figure}

We first present a bifurcation diagram in the $(K_c,k_{nt})$ plane in Fig. \ref{KckntBifurcation} (a), where $\tau=10$. In this example, the region of unique stable limit cycles (Region II) is L-shaped. Generally speaking, for large $k_{nt}$ values ($k_{nt}\geqslant 0.1$), $K_c$ must be small for oscillations to occur, and the part of Region II corresponding to these large $k_{nt}$ values is very narrow, forming the vertical part of an L-shape; for small $k_{nt}$ values ($k_{nt}< 0.1$), $K_c$ can be much larger, and the corresponding part of Region II forms the horizontal part of the L-shape. Note that the same trends of how $k_{nt}$ affects the region of oscillation may be seen in Fig. \ref{tauKcHopfCurves}.

We then present Hopf bifurcation curves of equilibria with various $\tau$ values in Fig. \ref{KckntBifurcation} (b), which depict a series of cross-sections of the Hopf bifurcation surface in the 3D parameter space parallel to the $(K_c,k_{nt})$ plane. We see that as $\tau$ increases, the vertical part of the L-shape decreases in height, meaning that oscillation is only possible for smaller values of $k_{nt}$. The area of the horizontal part of the L-shape also decreases as $\tau$ increases. In general, increasing $\tau$ causes the oscillation region to shrink in the $(K_c,k_{nt})$ plane. Note that Bautin bifurcations (marked by dashed curves) are also visible in Fig. \ref{KckntBifurcation}.

\begin{figure}[h]
 \includegraphics[width=0.8\textwidth]{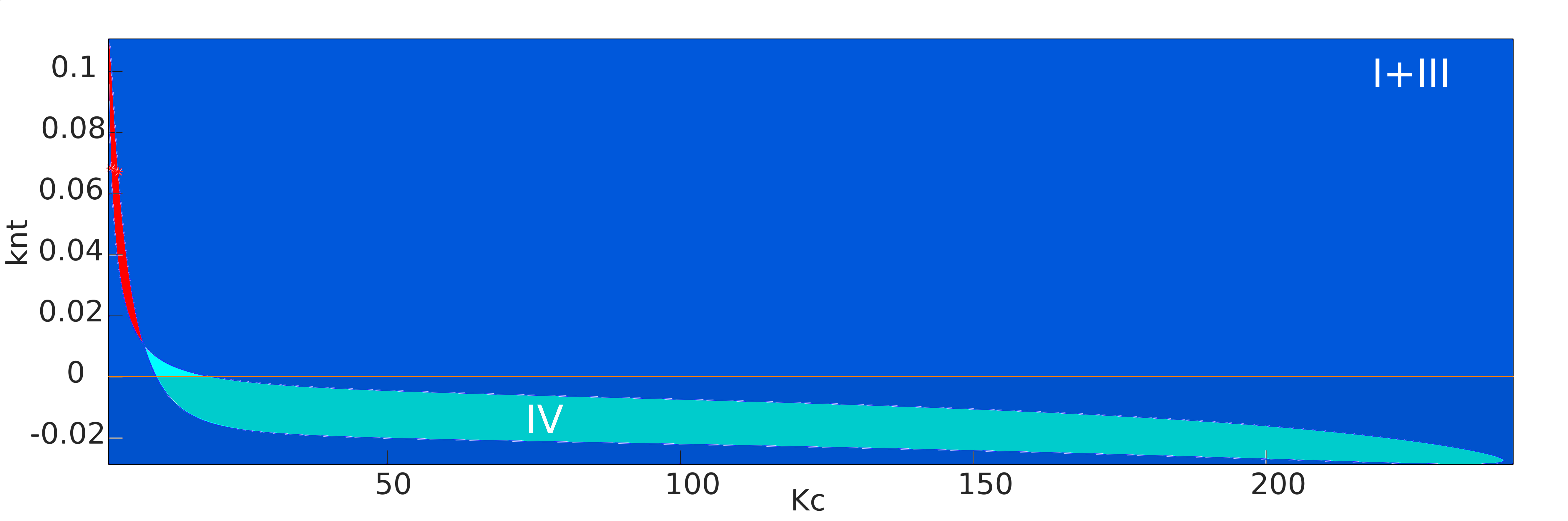}
 \caption{$(K_c,k_{nt})$ plane partitioned by Hopf bifurcation curve of equilibria, with $\tau=60$. Blue (Regions I and III), one stable equilibrium, except for very narrow regions along the Hopf bifurcation curve, where stable and unstable limit cycles coexist. Red (Region II), one stable limit cycle. Cyan (Region IV), two stable equilibria and one unstable equilibrium (bistability). The half-plane $k_{nt} < 0$ is shaded.}
 \label{KckntBistability}
\end{figure}

Although we have discussed the mechanism for the bistability of equilibria in Section \ref{BistabilityinEquilibria}, it is interesting to further study this bistability in the bifurcation analysis for $(K_c,k_{nt})$. Bistable regions may already be seen in Fig. \ref{KckntBifurcation}(a), but they are better visualized by allowing $k_{nt}$ to go beyond its biologically realistic range of nonnegative values. We see from Fig. \ref{KckntBistability} that, if the continuation of the Hopf bifurcation curve proceeds into the $k_{nt}<0$ half-plane, there is a relatively wide region of bistability (Region IV).

\subsection{Role of noise and robustness of oscillations}

The ODE model we studied in the previous section arises from a stochastic system of reactions \cite{chae2023spatially} that is influenced by finite-size noise, as well as various other noise sources. We investigate here the role of noise in the oscillations we reported, particularly near their onset, and the possible role played by the Bautin bifurcation in the regularity of the oscillations. Two types of noise may be relevant in this study: the chemical Langevin equation approximation \cite{anderson2011continuous,anderson2015stochastic}, or a simple additive noise model, which is considered hereafter. We now investigate the dynamics of the system of stochastic differential equations with additive noise:

\begin{subequations}\label{sdeexp2a}
 \begin{align}
  \label{sde1exp2a} dc_{no} &= [k_{vn}-f_p(c_{no},c_{nop},K_c,\tau)+f_{dp}(c_{no},c_{nop},\tau)-c_{no}A_{cyto}(k_{vn}+k_{vcy}-k_{nt}c_{nop})]dt \notag\\
  &+ \sigma dW_t^{no}, \\
  \label{sde2exp2a} dc_{nop} &= [f_p(c_{no},c_{nop},K_c,\tau)-f_{dp}(c_{no},c_{nop},\tau)-k_{nt}c_{nop}-c_{nop}A_{cyto}(k_{vn}+k_{vcy}-k_{nt}c_{nop})]dt \notag\\
  &+ \sigma dW_t^{nop},
 \end{align}
\end{subequations}
where $W_t^{no}$ and $W_t^{nop}$ are independent Brownian motions, and the parameter $\sigma$ controls a common noise level. 

\subsubsection{Stochastic resonance}

\begin{figure}[h]
 \includegraphics[width=\textwidth]{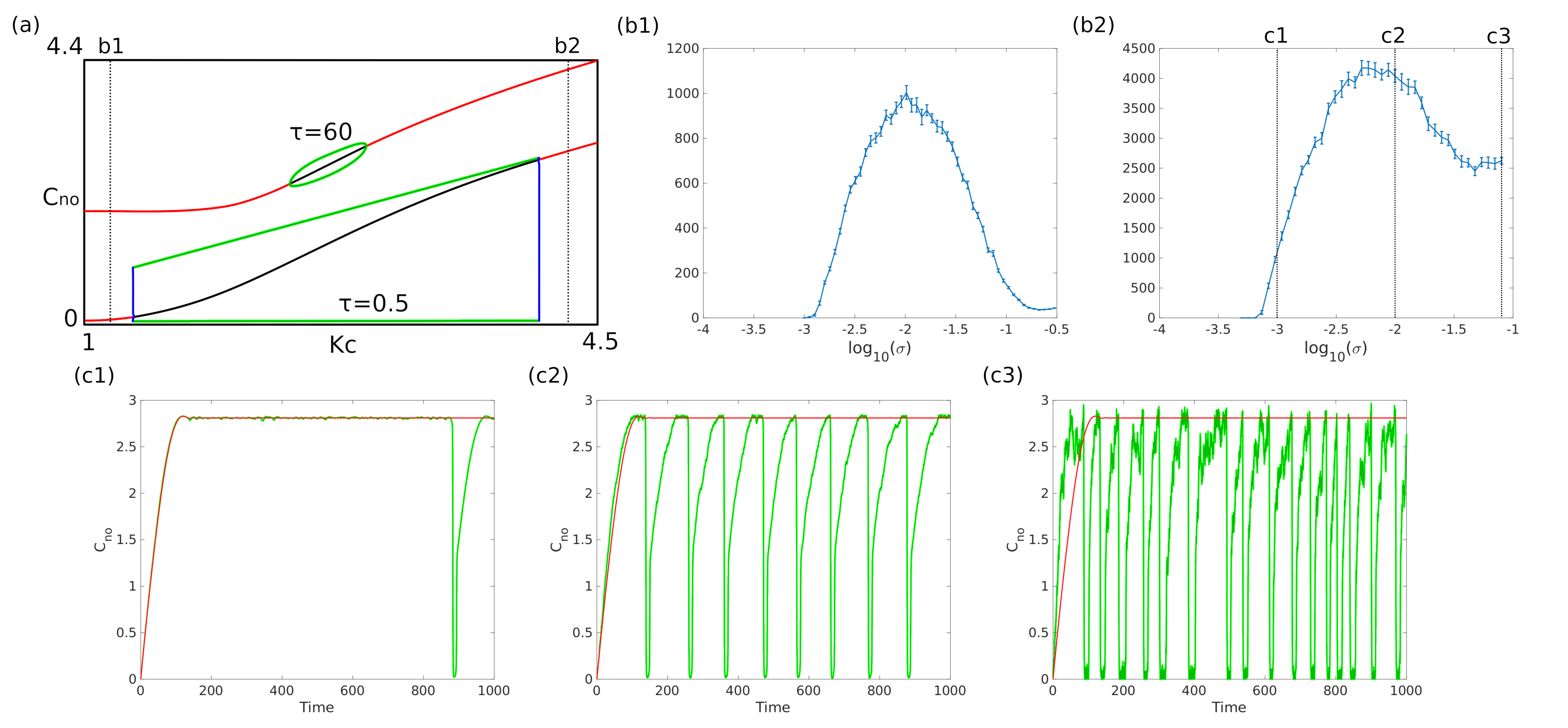}
 \caption{Stochastic resonance near Hopf bifurcation points. \textbf{(a)}, bifurcation with respect to $K_c$ at $\tau=0.5$ and $\tau=60$. Vertical lines label the $K_c$ values used in stochastic resonance examples. \textbf{(b1)(b2)}, maximal amplitude of the Fourier transform as a function of the noise level $\sigma$ at $K_c=1.25$ (onset of oscillations) and $K_c=4.2$ (end of oscillations), respectively. Vertical lines label the noise levels in sample trajectories. \textbf{(c1)(c2)(c3)}, deterministic (red) and stochastic (green) trajectories for $K_c=4.2$ and $\sigma=0.001$ (onset of stochastic resonance), $\sigma=0.01$ (peak of stochastic resonance), and $\sigma=0.08$.}
 \label{StochasticResonance}
\end{figure}

We first investigate the dynamics of the stochastic system for various noise levels $\sigma$, particularly near the onset of oscillations. To identify the possible presence of stochastic resonance \cite{muratov2005self} where noise enhances the reliability of oscillations in a system that would not oscillate without noise, we quantified the coherence of the responses of the stochastically forced system by estimating the maximal amplitude of the Fourier transform of the $c_{no}$ trajectory \cite{touboul2018complex} (Fig. \ref{StochasticResonance} (b)).

As shown in Fig. \ref{StochasticResonance} (a), stochastic resonance occurs when $k_c$ is slightly less than the first Hopf bifurcation point (onset of oscillations in the $K_c$ bifurcation diagram) and slightly greater than the second limit point of the Hopf bifurcation point (end of oscillations in the $K_c$ bifurcation diagram).

\begin{figure}[h]
 \includegraphics[width=0.5\textwidth]{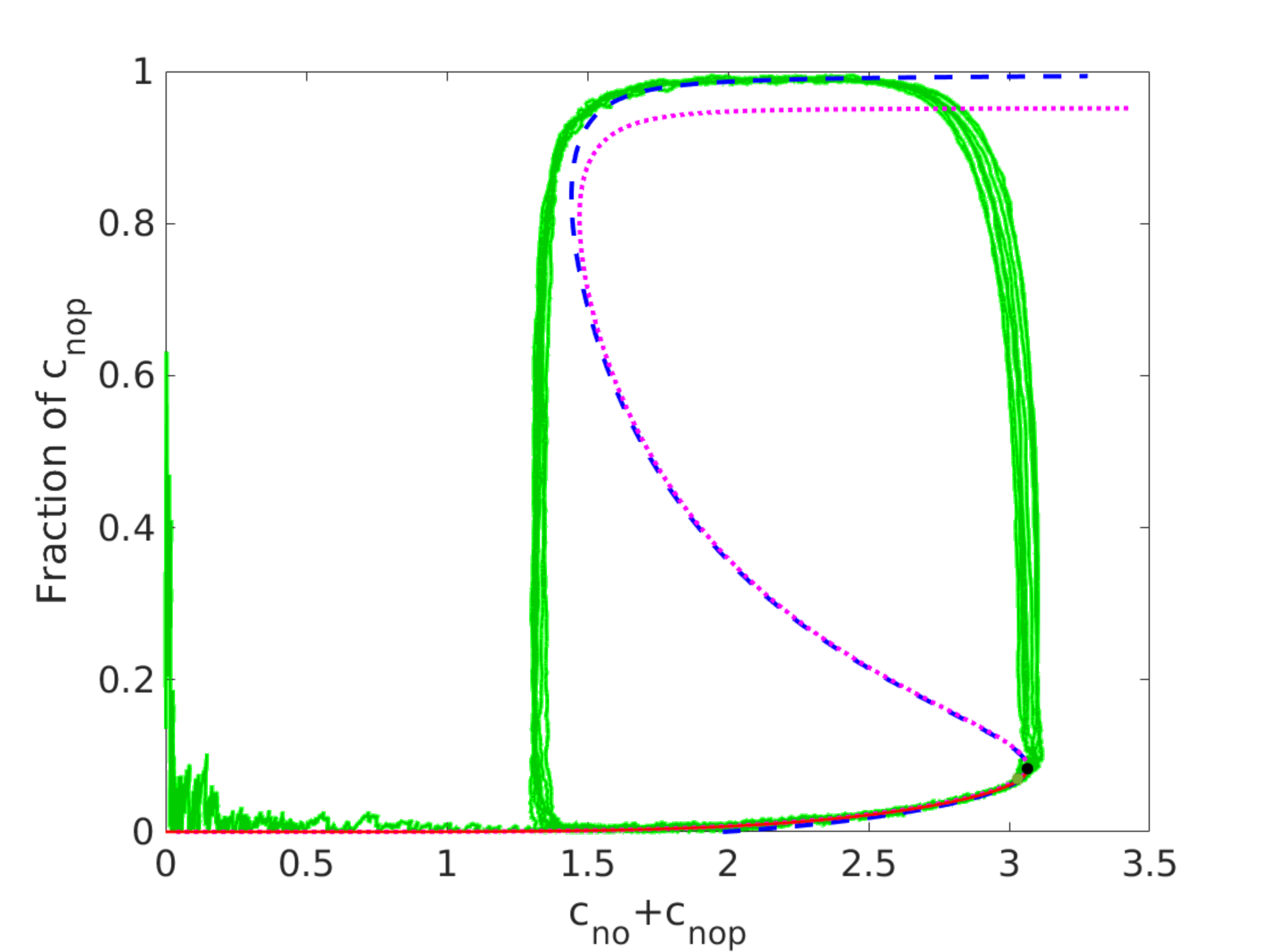}
 \caption{Phase plane illustrating the mechanism of stochastic resonance. Red solid line, deterministic trajectory; green solid line, stochastic trajectory; dashed line, $c_{no}$-nullcline; dotted line, $c_{nop}$-nullcline. $K_c=4.2$, $\tau=0.5$, $k_{nt}=0.1$, $\sigma=0.01$.}
 \label{SRMechanism}
\end{figure}

The mechanism generating stochastic resonance is the S-shape of the nullclines resulting from bistable phosphorylation, as noise may drive the trajectory past the equilibrium and beyond the fold points of the response curve, leading to cycles (Fig. \ref{SRMechanism}). When the noise level is low, the trajectory is highly unlikely to fluctuate beyond the equilibrium and cycles rarely occur (Fig. \ref{StochasticResonance} (c1)). At an intermediate level, noise induces oscillations with a regular period, similar to the vanished limit cycle (Fig. \ref{StochasticResonance} (c2)). If we further increase the noise level, the trajectory becomes noisy and the oscillation period becomes irregular (Fig. \ref{StochasticResonance} (c3)).

Fig. \ref{StochasticResonance} shows that compared to the corresponding deterministic system (\ref{odeexp2}), the stochastic system (\ref{sdeexp2a}) uses noise to generate oscillations across a wider range of parameter values. When the phosphorylation threshold $K_c$ is slightly outside the deterministic region of oscillation, the stochastic system can still oscillate and generate rhythms. In summary, bistable phosphorylation enhances the robustness of rhythm generation via stochastic resonance.

\subsubsection{Robustness of oscillation periods}

\begin{figure}[h]
 \includegraphics[width=0.9\textwidth]{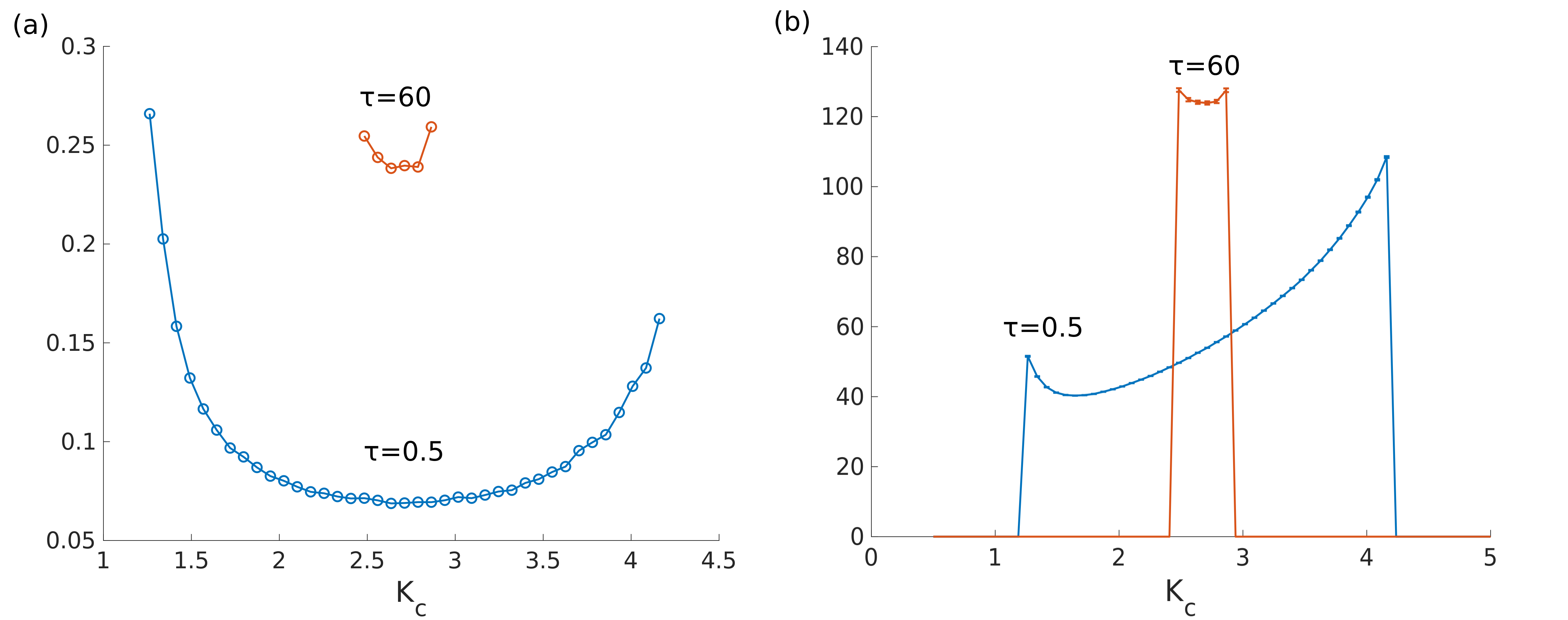}
 \caption{Robustness of oscillation periods. \textbf{(a)}, coefficient of variation of oscillation period as a function of $K_c$ at $\tau=0.5$ and $\tau=60$. \textbf{(b)}, corresponding average oscillation period as a function of $K_c$. Noise level $\sigma=0.01$.}
 \label{CV_period}
\end{figure}

To further study the robustness of oscillations under noise and examine whether the Bautin bifurcation described in Section \ref{tauKcBifurcationSec} enhances the robustness of oscillations, we calculate the coefficients of variation of oscillation periods as a function of the phosphorylation threshold $K_c$ at $\tau=0.5$ and $\tau=60$, each from 50 sample trajectories of the stochastic system (\ref{sdeexp2a}) at an intermediate noise level ($\sigma=0.01$).

For $\tau=0.5$, which lies before the Bautin bifurcation (Fig. \ref{tauKcBifurcation}), subcritical Hopf bifurcations lead to the sudden onset of oscillations with large amplitude (Fig. \ref{StochasticResonance}(a)). The coefficient of variation curve (Fig. \ref{CV_period}(a), blue) has a relatively flat section where the coefficient of variation of periods remains small ($<0.1$) for a wide range of $K_c$ values. This flat section indicates that such large-amplitude oscillations are robust against noise. The coefficient of variation increases significantly at $K_c$ values near the boundary of the deterministic region of oscillation, especially at $K_c$ values where stochastic resonance occurs. Overall, the coefficient of variation curve is U-shaped with a nearly constant flat section. However, the corresponding average period curve (Fig. \ref{CV_period}(b), blue) does not have a near-constant section, indicating that variation in the phosphorylation threshold $K_c$ may significantly affect the oscillation period.  

For  $\tau=60$, which lies after the Bautin bifurcation (Fig. \ref{tauKcBifurcation}), supercritical Hopf bifurcations generate oscillations with much smaller amplitude (Fig. \ref{StochasticResonance}(a)). The coefficient of variation (Fig. \ref{CV_period}(a), red) lies in a much higher interval ($\sim 0.25$), indicating that these small-amplitude oscillations are less robust against noise. In contrast, the corresponding average period (Fig. \ref{CV_period}(b), red) shows less variation than in the small $\tau$ case.

Overall, from the perspective of rhythm generation, small $\tau$ values and the presence of subcritical Hopf bifurcations enhance the robustness of oscillations against noise; however, the oscillation period still depends strongly on the parameter values.

\section{Discussion}

The bifurcation analysis described in Section \ref{results} shows that the bistable phosphorylation of nuclear cargo proteins may generate oscillations in protein concentrations. Our findings suggest that bistable phosphorylation may serve as a potential rhythm-generating mechanism that operates independently, without the need for other mechanisms such as negative feedback \cite{chae2023spatially,kim2014molecular,kim2021tick}. However, for each fixed value of $k_{nt}$ (Fig. \ref{tauKcBifurcation} and \ref{tauKcHopfCurves}) and $\tau$ (Fig. \ref{KckntBifurcation}), the area of the oscillation region is small compared to the 2D parameter plane explored. On the basis of these observations, bistable phosphorylation requires a somewhat specific set of parameters to generate oscillations. It is possible that the region of oscillation in the parameter space increases when the parameters are restricted to biologically realistic values, though parameter identification from experimental observations is beyond the scope of this paper.

Our study of the stochastic model shows that bistable phosphorylation may improve the robustness of rhythm generation against noise in two ways. One is by inducing stochastic resonance and using noise to extend the region of oscillation in parameter space; the other is by generating large-amplitude oscillations with more stable periods under noise via a subcritical Hopf bifurcation. Such findings supplement the study of robustness-enhancing features of bistable phosphorylation in circadian timekeeping \cite{chae2023spatially}. In the context of our model of nucleocytoplasmic transport and cell growth, this elucidates a potential mechanism for significant fluctuations in the nuclear-to-cytoplasmic volume ratio. Crucially, this bistable phosphorylation mechanism remains viable even in the biologically relevant regime of a large number of biomolecules.

\subsection{Limitations of the model and directions for future research}

\subsubsection{Phosphorylation of nuclear cargo}
We assumed that all cargo proteins must be phosphorylated prior to nuclear import and that the phosphorylation dynamics is modeled by a bistable response curve \cite{chae2023spatially}. These simplifying assumptions reduce the complexity of the model and the number of variables but may not be biologically realistic.

In reality, the role and dynamics of phosphorylation in nucleocytoplasmic transport are far more complex. Phosphorylation may both promote and inhibit the nuclear import of cargo proteins through numerous mechanisms \cite{nardozzi2010phosphorylation}, and may also regulate nuclear export \cite{jans1996regulation}. Given the complexity and diversity of the phosphorylation regulation pathways in nucleocytoplasmic transport, the dynamics of phosphorylation may be more complicated than the bistable response curve model suggests. Therefore, more realistic models must categorize nuclear cargo proteins by different types of phosphorylation dynamics and incorporate more details into the associated phosphorylation pathways. However, such models are likely to be significantly more complex and less amenable to the bifurcation analysis presented here.

\subsubsection{Interactions of bistable phosphorylation and negative feedback}
Our simplified model describes the dynamics of bistable phosphorylation in isolation. As stated above, the original circadian rhythm models \cite{chae2023spatially,kim2014molecular,kim2021tick} also include various negative feedback components as alternative mechanisms for rhythm generation. To our knowledge, while there are detailed bifurcation analyses of the Kim-Forger model with negative feedback \cite{pei2024three}, there is currently no bifurcation analysis of circadian rhythm models with bistable phosphorylation and negative feedback. We plan to further study the original model in \cite{chae2023spatially} from the perspective of dynamical systems and examine deterministic and stochastic dynamics. Such studies will explore the rich dynamics of the interactions between bistable phosphorylation and negative feedback.

\subsubsection{Theoretical aspects of bifurcation analysis}
We have performed all the bifurcation analyses in this paper using a combination of existing and freely available numerical software. While these numerical tools can compute regions in parameter space corresponding to different dynamical behaviors of the system, some of these calculations are rather complicated to perform in practice and require fine-tuning of algorithmic parameters. More importantly, numerical bifurcation analysis cannot fully identify the mechanism behind certain dynamical behaviors. For example, as seen in Section \ref{2parametertauKc}, it is challenging to identify the regions in which stable and unstable limit cycles coexist. We also see from this example that, although numerical bifurcation analysis shows that these regions originate from Bautin bifurcations near $\tau=40$ in Fig. \ref{tauKcHopfCurves}, it does not provide a straightforward explanation for why these coexistence regions become narrower and finally disappear as $\tau\to 0$.

One potential direction for future studies is to characterize the dynamics through the theory of geometric bifurcations \cite{barrio2024exploring}, which uses Morse theory to classify the 2D bifurcation surfaces in 3D parameter spaces. We expect that the theory of geometric bifurcations, together with other theoretical tools for bifurcation analysis, could provide further knowledge of the global topology and geometry of the bifurcation surface in the $(\tau,K_c,k_{nt})$ space and provide a more systematic understanding of the various dynamical behaviors of the model.

\section*{Acknowledgments}

We acknowledge helpful discussions with Jae Kyoung Kim. We acknowledge funding from NSF grant MCB-2213583 to XB and TGF.

\printbibliography

\end{document}